\documentclass[a4paper, 12pt]{amsart}

\def\nnewpage{\newpage}

\usepackage[utf8]{inputenc}
\usepackage[american]{babel}

\usepackage[draft=true]{minted}

\usepackage{amsmath,amssymb,amsthm}
\usepackage{bbm}    

\usepackage{bm}
\usepackage{accents}
\usepackage{mathrsfs}
\usepackage{mathtools}
\usepackage{fixmath}
\usepackage{paralist}
\usepackage{enumerate}
\usepackage{xspace}

\usepackage{stmaryrd}

\usepackage{booktabs}
\usepackage{array}
\usepackage{graphicx}
\usepackage{pgf,tikz}
\usepackage{array}
\usepackage{hyperref}
\usepackage{cleveref}
\usepackage{enumitem}
\usepackage{csquotes}

\usepackage{pgf,tikz}
\usetikzlibrary{arrows}
\usetikzlibrary{patterns}
\usetikzlibrary[positioning]
\usetikzlibrary[fit]
\usetikzlibrary{shapes.geometric}
\usetikzlibrary{calc}

\usepackage{pgfplots}
\pgfplotsset{compat=1.17}

\newcommand\CC{{\mathbb C}}

\newcommand\NN{{\mathbb N}}

\newcommand\QQ{{\mathbb Q}}
\newcommand\RR{{\mathbb R}}

\newcommand\ZZ{{\mathbb Z}}

\newcommand\1{{\mathbf 1}}

\newcommand\smallSetOf[2]{\{#1 \mid #2\}}

\DeclareMathOperator\conv{conv}

\DeclareMathOperator\Gal{Gal}
\DeclareMathOperator\DP{DP}

\newtheorem{theorem}{Theorem}

\theoremstyle{remark}

\newtheorem{example}[theorem]{Example}

\newtheorem{remark}[theorem]{Remark}
\newtheorem{question}[theorem]{Question}

\newcommand\Julia{\texttt{Julia}\xspace}
\newcommand\polymake{\texttt{polymake}\xspace}
\newcommand\Singular{\texttt{Singular}\xspace}
\newcommand\ANTIC{\texttt{ANTIC}\xspace}
\newcommand\GAP{\texttt{GAP}\xspace}

\newcommand\OSCAR{\texttt{OSCAR}\xspace}
\newcommand\Pari{\texttt{Pari}\xspace}
\newcommand\SageMath{\texttt{SageMath}\xspace}
\newcommand\Magma{\texttt{Magma}\xspace}

\newcommand\latte{\texttt{LattE}\xspace}
\newcommand\normaliz{\texttt{Normaliz}\xspace}
\newcommand\polydb{\texttt{polyDB}\xspace}

 \usepackage[colorinlistoftodos,bordercolor=orange,backgroundcolor=orange!20,linecolor=orange,textsize=scriptsize]{todonotes}
\usepackage{fullpage}

\begin{document}

\title[Galois groups of Ehrhart polynomials in \OSCAR]{Computing Galois groups of Ehrhart \\ polynomials in \OSCAR}

\author{Claus Fieker}
\address{Technische Universität Kaiserslautern}
\email{fieker@mathematik.uni-kl.de}

\author{Tommy Hofmann}
\address{Universität Siegen}
\email{tommy.hofmann@uni-siegen.de}

\author{Michael Joswig}
\address{Technische Universität Berlin, Chair of Discrete Mathematics/Geometry, Max-Planck Institute for Mathematics in the Sciences, Leipzig}
\email{joswig@math.tu-berlin.de}

\begin{abstract}
  We report on an implementation of Galois groups in the new computer algebra system \texttt{OSCAR}.
  As an application we compute Galois groups of Ehrhart polynomials of lattice polytopes.
\end{abstract}


\maketitle

\section{Introduction}
\OSCAR is an acronym for \enquote{Open Source Computer Algebra Resource} \cite{OSCAR,OSCAR-book}.
This is a new computer algebra system, written in \Julia, which combines and extends the full feature set of \ANTIC \cite{nemo}, \GAP \cite{GAP4}, \polymake \cite{DMV:polymake} and \Singular \cite{Singular}. %
The implementation is an ongoing major collaborative effort lead by the transregional collaborative research center (SFB-TRR) 195 \enquote{Symbolic Tools in Mathematics and their Application}, which is funded by the German Research Foundation (DFG).
At the current version 0.8.2 the system is still in its infancy; yet it offers a wide range of functions already.
As a showcase here we give an account of computations in \OSCAR which require methods from algebra, number theory, group theory and polyhedral geometry.
\OSCAR is registered as a \Julia package, which makes the installation of a released version via \Julia's package manager trivial.
Additionally, the entire source code and its ongoing development can be followed at \url{https://github.com/oscar-system}.
This also includes documentation, unit tests and continuous integration.
The specific computations shown here are documented at \url{https://github.com/micjoswig/oscar-notebooks}.

A \emph{lattice polytope}, $P$, is the convex hull of finitely many points with integral coordinates.
The function $L_P(t)$, for $t\in\NN$, counts the lattice points in the dilations $t\cdot P$.
Ehrhart's Theorem says that $L_P(t)$ is a univariate polynomial of degree $d=\dim P$ in the parameter $t$; cf.\ \cite[\S 3.4]{Beck+Robins:2015}. %
Research on Ehrhart polynomials abounds for their connections to combinatorics, algebraic geometry and beyond.
Here we use \OSCAR to study the Galois groups of Ehrhart polynomials.
To the best of our knowledge this is a new idea.

\section{Computing Galois Groups}

The computation of Galois groups of rational polynomials is among the fundamental problems of algorithmic algebraic number theory
as formulated by Zassenhaus in \cite{Zassenhaus:1982}.
In the following we will give a brief overview about the standard algorithm for solving this task as well as the state of current implementations.

We will first explain what is meant by computing the Galois group of a rational polynomial.
To this end let $f \in \QQ[x]$ be a polynomial of degree $n$.
Recall that the Galois group $\Gal(f)$ of $f$ is defined as the group of automorphisms of $\Gal(N/\QQ)$, where $N$ denotes a splitting field of $f$,
i.e., a minimal field such that $f$ decomposes over $N$ into linear factors.
Denote by $S = \{\alpha_1,\dotsc,\alpha_n\}$ the set of roots of $f$ in any field containing $\QQ$.
After choosing an ordering on the roots $S$, via its action on $S$ the Galois group $\Gal(f)$
gives rise to a faithful permutation presentation $G \to
S_n$.
The equivalence class of this presentation is independent of the chosen set $S$ or its ordering.
Thus when speaking of computing $\Gal(f)$, we mean the computation of a permutation group $G \subseteq S_n$ such that
the action of $\Gal(f)$ on a set of roots in any splitting field is permutation equivalent to the natural action of $G$.
Note that this implies $\Gal(f) \cong G$ as groups, but permutation equivalence is in general stronger.

We now sketch an algorithm of Stauduhar~\cite{Stauduhar:1973} and its improvements by Fieker and Klüners~\cite{Fieker+Kluners:2014}
for finding the Galois group $\Gal(f) \leq S_n$, where $f \in \QQ[x]$ has degree $n$.

\begin{enumerate}
  \item
    Determine approximations $\hat S = \{\hat \alpha_1, \dotsc, \hat \alpha_n\} \subseteq L$ of the roots of $f$ in a suitable field extension $\QQ \subseteq L$. One may take $L = \CC$ and employ well-known methods to approximate roots of univariate polynomials, but it is also possible to use algebraic closures of $p$-adic fields $\QQ_p$.
  \item
    Assume that we know a subgroup $G \leq S_n$ for which the inclusion $\Gal(f) \leq G$ holds (at the beginning $G = S_n$). To test whether $\Gal(f)$ is equal to $G$, it is sufficient to test whether $\Gal(f) \leq H$ for any of the finite maximal subgroups $H \leq G$.
    Thus let $H \leq G$ be a subgroup. Denote by $F \in \ZZ[x_1,\dotsc,x_n]$ a $G$-relative $H$-invariant polynomial, that is, $F$ equals $F^\sigma := F(x_{\sigma(1)},\dotsc,x_{\sigma(n)})$ and the stabilizer of $F$ in $G$ equals $H$.
    Now after certain technical properties are checked, the main result of Stauduhar states that $\Gal(f) \leq H$ if and only if $F(\alpha_1,\dotsc,\alpha_n)$ is an element of $\ZZ$.
    As this latter task can performed by evaluating $F^\sigma(\hat \alpha_1,\dotsc,\hat \alpha_n)$ to a high enough precision, we can decide whether $\Gal(f) \leq H$ rigorously.
    By performing this task iteratively, the algorithm traverses down the lattice of subgroups of $S_n$ until the Galois group $\Gal(f) \leq S_n$ is determined.
\end{enumerate}

\begin{remark}
  While the algorithm described above yields a procedure to compute Galois groups that works well in practice for input of small size,
  we want to point out several improvements, which are necessary to render this algorithm practical.
  \begin{enumerate}
    \item
      While we formulated the algorithm for arbitrary polynomials, it is advisable to first determine a factorization of $f$ into irreducible polynomials and compute the Galois group for each of these polynomials first. This yields
      two advantages:
      For irreducible polynomials, the Galois is automatically transitive, hence fewer 
      candidates have to be checked and secondly, for pairs of transitive
      groups more sophisticated methods for finding invariants are available
      than in the generic case.
    \item
      Instead of testing $\Gal(f) \leq H$ for all maximal subgroups $H$ of $G$, one can restrict to representatives of conjugacy classes of maximal subgroups of $G$.
      Given such a representative $H$, one then tests for a set of coset representatives $\tau$ of $H$ in $G$ whether $F^\tau(\alpha_1,\dotsc,\alpha_n)$ is an integer, which is equivalent to $G \leq H^\tau$.
      Determining the conjugacy classes of maximal subgroups and coset representatives of permutation groups are well-studied problems in algorithmic group theory, with efficient solutions.
    \item
      One bottleneck of the algorithm is finding the $G$-relative $H$-invariant polynomial $F$. Generic methods described in \cite{Fieker+Kluners:2014} yield a polynomial whose number of terms grows linearly in the index $[G : H]$, making this step quite costly for groups with maximal subgroups of large index.
      More advanced techniques have been introduced by Elsenhans~\cite{Elsenhans:2017}.
    \item To deal with large index subgroups, the idea of short cosets has been introduced, where additional information is used to reduce the set of coset representatives that need to be considered. For example, if $f$ has a pair of complex conjugate roots, the complex
      conjugation has to be an element of the Galois group of $f$, and hence only
      $\sigma$ such that $H^\sigma$ contains the complex conjugation need to be
      tested.
    \item To make use of highly optimized algorithms for finding subfields of number fields, instead of starting with the full symmetric
      group, one can start with the intersection of suitable wreath products:
      If $f$ is irreducible, each subfield of the stem field $\QQ[x]/(f)$ of $f$ yields a block system for the
      unknown Galois group, hence a maximal subgroup supporting this.
  \end{enumerate}
\end{remark}

An implementation of the algorithm, which performs quite well in practice, is available in \OSCAR (version 0.8.2).
The implementation exploits the tight integration between the group theoretic component built upon \GAP and the number theoretic functionality provided by \ANTIC.
Section~\ref{sec:Ehrhart} below comprises concrete computations in the context of Ehrhart polynomials.

We end with comparisons to other computer algebra packages with the functionality to compute Galois groups of (a restricted set of) polynomials in $\QQ[x]$.
As is apparent from the description of the algorithm, critical steps require the use of sophisticated techniques from both algorithmic group and number theory.
This requirement already narrows down possible computer algebra software able to handle this scenario.

\begin{enumerate}
  \item
    Computing Galois groups of rational polynomials has been available in the closed source system \Magma~\cite{Magma} for quite some time.
    In comparison to \OSCAR, due to various optimizations, for increasing parameters \Magma will outperform the current implementation in \OSCAR.
    Yet, for parameter ranges of practical relevance, the performance profiles are similar.
  \item
    Based on the work of Eichenlaub~\cite{Eichenlaub:1996}, the number theory package \Pari (version 2.13.1) provides an implementation to determine Galois groups for polynomials $f$ with degree $\deg(f) \leq 11$, which in addition need to be irreducible.
    Neither assumption is met in our application: Ehrhart polynomials may be reducible, and they may have arbitrarily high degrees.
  \item
    Implementations of algorithms of Soicher--McKay \cite{Soicher+McKay:1985} and Hulpke~\cite{Hulpke:1999} in \GAP allow to determine Galois groups of irreducible polynomials of degree up to 15.
  \item
    The computer algebra system \SageMath~\cite{sagemath} (version 9.5) includes interfaces to \Pari and \GAP, which are subject to the same restrictions ($f$ needs to be irreducible with $\deg(f) \leq 11$ or $\deg(f) \leq 15$, depending on which interface is used).

\end{enumerate}

Thus, for the first time, Galois groups of arbitrary rational polynomials can be determined using a freely available open source computer algebra package.
For example, the Galois group of the Ehrhart polynomial of the 16-dimensional Fano simplex of order $2^8 \cdot 8!$ (see Section~\ref{sec:Ehrhart}) can be determined using \OSCAR within 10 minutes, a computation which was previously impossible without \Magma.

\section{Ehrhart Polynomials of Lattice Polytopes}\label{sec:Ehrhart}

Let $P\in\RR^d$ be a $d$-dimensional lattice polytope.
Then the function $L_P(t)=|tP\cap\ZZ^d|$ is a polynomial in $t$ of degree $d$; this is called the \emph{Ehrhart polynomial} of $P$.
Numerous facts are known about Ehrhart polynomials: e.g., the leading coefficient is the volume of $P$, the constant term equals one, and the polynomial $d! \cdot L_P(t)$ has integral coefficients.
Yet many questions remain open.
For instance, it is unclear which polynomials arise as Ehrhart polynomials \cite[Open Problem 3.43]{Beck+Robins:2015}.
We refer to the monographs of Beck and Robins \cite{Beck+Robins:2015} and Ewald \cite{Ewald:1996} for the general background.

Computing Ehrhart polynomials is a subtle task, which may involve triangulating the input and computing lattice points; cf.\ \cite[\S3.4]{Beck+Robins:2015} and \cite{Barvinok:2008}.
Note that even computing the leading coefficient, which is the volume, is $\#P$-hard \cite{Dyer+Frieze:1988}.
Standard implementations include \latte \cite{latte} and \normaliz \cite{normaliz}; both are available through interfaces in \polymake and thus also in \OSCAR.
For comparisons and further implementations see \cite{polymake:2017}.

\paragraph*{Galois groups.}
Here we propose to investigate the Galois groups of Ehrhart polynomials, considered as polynomials over some algebraically closed extension of $\QQ$, such as $\CC$ or the algebraic closure of the $p$-adics.
This extends ongoing research on the zeros of Ehrhart polynomials; see, e.g., \cite{BDDPS:2005,Bey+Henk+Wills:2007} and their references as well as \cite[Open Problem 3.44]{Beck+Robins:2015}.
We start with a first example, which is straightforward.
Its purpose is to exhibit that Ehrhart polynomials may be reducible.
\begin{example}\label{exmp:cube}
  Let $C$ be the $d$-dimensional cube $[-1,1]^d$.
  Its Ehrhart polynomial reads $L_C(t)=(2t+1)^d$, whose Galois group is trivial.
\end{example}

The next example is more interesting.
\begin{example}\label{exmp:fano-simplex}
  Consider the $d$-simplex $S=\conv\{e_1,e_2,\dots,e_d,-e_1-e_2-\dots-e_d\}$, whose normalized volume equals $d+1$.
  Here the Ehrhart polynomial is
  \[
    L_{S}(t) \ = \ \binom{t+d+1}{d+1} - \binom{t}{d+1} \enspace .
  \]
  Its roots are determined in \cite[Theorem 1.7]{Bey+Henk+Wills:2007}.
\end{example}
The polytope $S$ in the example above is sometimes called the \emph{Fano simplex}.
Employing the results of \cite{Bey+Henk+Wills:2007} one can prove the following.
\begin{theorem}
  The Galois group $\Gal(L_{S}(t))$ of the $d$-dimensional Fano simplex $S$ is isomorphic to the wreath product $C_2 \wr S_{k}$, where $k=\lfloor d/2 \rfloor$.
  Its order is $2^k k!$.
\end{theorem}

We show how to compute the Ehrhart polynomial and its Galois group for $d=14$.
This requires an installation of \Julia (version $\geq 1.6$), with \OSCAR installed via \mintinline{jl}{using Pkg; Pkg.add("Oscar")}.
Timings are taken on an Apple M1, 16 GB and macOS version 11.6.
\begin{minted}{jlcon}
julia> using Oscar;

julia> S = fano_simplex(14)
A polyhedron in ambient dimension 14

julia> @time G = galois_group(ehrhart_polynomial(S))
 0.840115 seconds (3.59 M allocations: 280.822 MiB, 6.75% gc time)
(<permutation group of size 645120 with 9 generators>, Galois Context for 15*x^14 + 105*x^13 + 24115*x^12 + 143325*x^11 + 7724717*x^10 + 37312275*x^9 + 725938785*x^8 + 2681453775*x^7 + 21964438496*x^6 + 56663366760*x^5 + 201186840400*x^4 + 310989260400*x^3 + 429952217472*x^2 + 283465647360*x + 87178291200 and prime 59)

julia> describe(G[1])
"C2 x ((C2 x C2 x C2 x C2 x C2 x C2) : S7)"
\end{minted}

\paragraph{Smooth Fano polytopes.}
We conclude this note with an analysis of a particularly interesting class of lattice polytopes.
Let $P\subset\RR^d$ be a $d$-dimensional lattice polytope with the origin in the interior.
Its \emph{polar} is $P^\circ = \smallSetOf{y\in\RR^d}{\langle x,y\rangle \leq 1 \text{ for all } x\in P}$, which is again a $d$-dimensional polytope.
The lattice polytope $P$ is \emph{reflexive} if $P^\circ$ is a lattice polytope \cite[\S2.4]{toric+varieties}.
Further, the lattice polytope is said to be \emph{smooth} if the primitive facet normal vectors in each normal cone form a $\ZZ$-basis of the ambient lattice $\ZZ^d$; cf.\ \cite[\S2.4]{toric+varieties}.
In particular, smooth polytopes are necessarily simple.
A \emph{smooth Fano polytope} is a smooth and reflexive lattice polytope.
Examples include the cube (Example~\ref{exmp:cube}) and the polar of the Fano simplex (Example~\ref{exmp:fano-simplex}).
Notice that sometimes in the literature smoothness is associated with the face fan of $P$ rather than the normal fan; see, e.g., \cite{Hegedus+Kasprzyk:2011,AssarfJoswigPaffenholz:2014,Hegedus+Higashitani+Kasprzyk:2019}.
In that case our definition applies to the polar $P^\circ$.
Interest in reflexive polytopes comes from toric algebraic geometry and its applications to mathematical physics.

\begin{example}\label{exmp:delpezzo}
  The reflexive polytope $\DP(d)=\conv\{\pm e_1,\pm e_2,\dots,\pm e_d,\pm \1\}$ is the \emph{del Pezzo polytope} of dimension $d$.
  Its polar $\DP(d)^\circ$ is smooth Fano if and only if $d$ is even.
\end{example}

The Ehrhart polynomials of the polars of the smooth Fano polytopes and their zeros have been investigated a lot; cf. \cite{Bey+Henk+Wills:2007,Hegedus+Kasprzyk:2011,Hegedus+Higashitani+Kasprzyk:2019} and their references.
Hegedüs and Kasprzyk \cite[Theorem 1.5]{Hegedus+Kasprzyk:2011} proved that for $P$ a smooth Fano polytope of dimension $d\leq 5$ the roots of the Ehrhart polynomial $L_{P^\circ}(t)$ of the polar lie on the \enquote{critical line} of those points in $\CC$ whose real part equals $-\tfrac{1}{2}$.
Moreover, Hegedüs and Kasprzyk also showed that, up to lattice equivalence, there are precisely four smooth Fano $6$-polytopes $P$ such that the roots of $L_{P^\circ}(t)$ do \emph{not} lie on the critical line.
In \cite[Example 1.7]{Hegedus+Kasprzyk:2011} these four polytopes are identified via the IDs 1895, 1930, 4853, and 5817 in the Graded Ring Database \cite{Obro:0704.0049,GRDB}, which stores the polars of all smooth Fano polytopes of dimension $\leq 6$.
The second but last example, which is the free sum $\DP(2)\oplus\DP(4)$, is special; cf.\ \cite[Theorem 7(iii)]{AssarfJoswigPaffenholz:2014}.
An \OSCAR computation reveals that Galois group of its Ehrhart polynomial is isomorphic to $S_3$.
In fact, $\DP(2)\oplus\DP(4)$ is the only polar of a smooth Fano $6$-polytope such that the Ehrhart polynomial has Galois group of order six.

\OSCAR provides direct access to the database \polydb \cite{polydb} which extends the classification from \cite{Obro:0704.0049} to dimension $\leq 8$.
The respective \polydb IDs of (the polars of) the four polytopes from \cite[Example 1.7]{Hegedus+Kasprzyk:2011} read: F.6D.0803, F.6D.0720, F.6D.3154, and F.6D.2616.
We can query \polydb to obtain a tally of the orders of the Galois groups of the Ehrhart polynomials of the polars of all 7622 smooth Fano $6$-polytopes:
\begin{minted}{jlcon}
julia> db = Polymake.Polydb.get_db();

julia> collection = db["Polytopes.Lattice.SmoothReflexive"];

julia> query = Dict("DIM"=>6);

julia> res = Polymake.Polydb.find(collection, query);

julia> E6 = [ ehrhart_polynomial(polarize(Polyhedron(P))) for P in res ];

julia> sort(MSet(order(galois_group(e)[1]) for e = E6).Dict)
OrderedCollections.OrderedDict{Int64, Int64} with 6 entries:
  4  => 623
  6  => 1
  8  => 22
  12 => 44
  16 => 310
  48 => 6622
  
\end{minted}

\nnewpage

\section{Conclusion}
The following question is immediate:
\begin{question}
  Which groups arise as Galois groups of Ehrhart polynomials?
\end{question}
It is a deep open question whether each finite group occurs as some Galois group of a rational polynomial; cf.\ \cite{Malle+Matzat:1999}.
While it seems to be difficult to characterize all Ehrhart polynomials it may still be the case that the above question is simpler.
For instance, it is known that each finite cyclic group occurs as the Galois group of some rational polynomial.
Can these be chosen as Ehrhart polynomials?

\section{Acknowledgments}

We are indebted to all contributors of \OSCAR.
Furthermore, we wish to thank Christian Haase, Gunter Malle and Benjamin Nill for useful conversations and remarks.
C.~Fieker and M.~Joswig have been supported by \enquote{Symbolic Tools in Mathematics and their Application} (TRR 195, project-ID 286237555).

\bibliographystyle{plain}
\bibliography{galois.bib}

\end{document}